\documentclass[10pt]{article}  
\usepackage{ijcas_packages}
\usepackage{fancyhdr}
\usepackage{bibentry}
\allowdisplaybreaks[1]

\newcommand{\squeezeup}{\vspace{-2.5mm}}

\newtheorem{prop}{Proposition}

\title{\LARGE \textbf{Constrained Geometric Attitude Control on \( \SO \)} }
\author{Shankar Kulumani* and Taeyoung Lee}%

\begin{document}
\maketitle

\begin{abstract}
This paper presents a new geometric adaptive control system with state inequality constraints for the attitude dynamics of a rigid body. 
The control system is designed such that the desired attitude is asymptotically stabilized, while the controlled attitude trajectory avoids undesired regions defined by an inequality constraint. 
In addition, we develop an adaptive update law that enables attitude stabilization in the presence of unknown disturbances. 
The attitude dynamics and the proposed control systems are developed on the special orthogonal group such that singularities and ambiguities of other attitude parameterizations, such as Euler angles and quaternions are completely avoided. 
The effectiveness of the proposed control system is demonstrated through numerical simulations and experimental results.
\end{abstract}

\section{Introduction}\label{sec:intro}
Rigid body attitude control is an important problem for aerospace vehicles, ground and underwater vehicles, as well as robotic systems~\cite{hughes2004,wertz1978}.
One distinctive feature of the attitude dynamics of rigid bodies is that it evolves on a nonlinear manifold.
The three-dimensional special orthogonal group, or \( \SO \), is the set of \( 3 \times 3 \) orthogonal matrices whose determinant is one.
This configuration space is non-Euclidean and yields unique stability properties which are not observable on a linear space.
For example, it is impossible to achieve global attitude stabilization using continuous time-invariant feedback~\cite{bhat2000}.



Attitude control is typically studied using a variety of attitude parameterizations, such as Euler angles or quaternions~\cite{shuster1993}.
Attitude parameterizations fail to represent the nonlinear configuration space both globally and uniquely~\cite{chaturvedi2011a}.
For example, minimal attitude representations, such as Euler angle sequences or modified Rodriguez parameters, suffer from singularities.
These attitude representations are not suitable for large angular slews.
In order to avoid singularities, the designer must carefully switch the chosen Euler angle sequence based on the operating region.
Another option is to artificially limit the operating region of the rigid body.
This ensures the system operates in a region free from singularities but limits the performance capabilities and ability to perform arbitrarily large angular maneuvers.
Quaternions do not have singularities but they double cover the special orthogonal group.
As a result, any physical attitude is represented by a pair of antipodal quaternions on the three-sphere.
An immediate implication of this ambiguity is that closed-loop stability properties derived using quaternions may not hold for the physical rigid body evolving on the true configuration space, namely the special orthogonal group.
During implementation, the designer must carefully resolve this non-unique representation in quaternion based attitude control systems to avoid undesirable unwinding behavior~\cite{bhat2000}.
This behavior is characterized by situations where the rigid body starts close to the desired attitude, yet the system unnecessarily rotates through a large angle in spite of a small initial error. 

Many physical rigid body systems must perform large angular slews in the presence of state constraints.
For example, autonomous spacecraft or aerial systems are typically equipped with sensitive optical payloads, such as infrared or interferometric sensors.
These systems require retargeting while avoiding direct exposure to sunlight or other bright objects.
In addition, many ground based attitude testing environments, such as air bearing platforms, must operate in the presence of physical obstructions.
Determining a satisfactory attitude control maneuver in the presence of state constraints is a challenging task.
The removal of constrained regions from the rotational configuration space results in a \textit{nonconvex} region.
The attitude control problem in the feasible configuration space has been extensively studied~\cite{bullo2004,MayTeePaCC11,LEEITAC15}.
However, the attitude control problem in the presence of constraints has received much less attention.

Several approaches have been developed to treat the attitude control problem in the presence of constraints.
A conceptually straightforward approach is used in~\cite{hablani1999} to determine feasible attitude trajectories prior to implementation.
The algorithm determines an intermediate point such that an unconstrained maneuver can be calculated for each subsegment.
Typically, an optimal or easily implementable on-board control scheme for attitude maneuvers is applied to maneuver the vehicle along these segments.
In this manner, it is possible to solve the constrained attitude control problem by linking several intermediary unconstrained maneuvers.
While this method is conceptually simple, it is difficult to generalize for an arbitrary number of constraints.
In addition, this approach is only applicable to problems where the selection of intermediate points are computationally feasible.

The approach in~\cite{frazzoli2001} involves the use of randomized motion planning algorithms to solve the constrained attitude control problem.
A graph is generated consisting of vertices from an initial attitude to a desired attitude. 
A random iterative search is conducted to determine a path through a directed graph such that a given cost, which parameterizes the path cost, is minimized.
The random search approach can only stochastically guarantee attitude convergence as it can be shown that as the number of vertices in the graph grow, the probability of nonconvergence goes to zero.
However, the computational demand grows as the size of the graph is increased and a new graph is required when constraints are modified. 
As a result, random search approaches are ill-suited to on-board implementation or in scenarios that require agile maneuvers.

Model predictive control for spacecraft attitude dynamics is another popular approach and has been studied in~\cite{guiggiani2014,kalabic2014,gupta2015}.
These methods rely on linear or non-linear state dynamics to repeatedly solve a finite-time optimal control problem.
Also known as receding horizon control, the optimal control formulation allows for a straight forward method to incorporate state and control constraints.
Computing the optimal control strategy over a moving horizon allows for a form of feedback type control rather than the more typical open-loop optimal control solution.
Due to the iterative nature of solving optimization problems, model predictive control methods are computational expensive and frequently apply direct optimization methods to solve the necessary conditions for optimality.
Therefore, these methods are complicated to implement and may not be suitable for real-time control applications.
  
Artificial potential functions are commonly used to handle kinematic constraints for a wide range of problems in robotics~\cite{rimon1992}.
The goal is the design of attractive and repulsive terms which drive the system toward or away from a certain obstacle, respectively.
The attractive function is designed to drive the system towards the desired state.
Similarly, a repulsive function is constructed such that the system is directed away from the constraints. 
The superposition of these functions allows one to apply standard feedback control schemes for stabilization and tracking.
More specifically, artificial potential functions have previously been applied to the spacecraft attitude control problem in~\cite{lee2011b,mcinnes1994}.
However, both of these approaches were developed using attitude parameterizations, namely Euler angles and quaternions, and as such, they are limited by the singularities of minimal representations or the ambiguity of quaternions.

This paper is focused on developing an adaptive attitude control scheme in the presence of attitude inequality constraints on \(\SO\).
We apply a potential function based approach developed directly on the nonlinear manifold \(\SO\). 
By characterizing the attitude both globally and uniquely on \(\SO\), our approach avoids the issues of attitude parameterizations, such as kinematic singularities and ambiguities, and is geometrically exact. 
A configuration error function on \(\SO\) with a logarithmic barrier function is proposed to avoid constrained regions. 
Instead of calculating a priori trajectories, as in the randomized approaches, our approach results in a closed-loop attitude control system. 
This makes it ideal for on-board implementation on UAV or spacecraft systems. 
In addition, unlike previous approaches our control system can handle an arbitrary number of constrained regions without modification.
This approach results in a conceptually simple obstacle avoidance scheme which extends the previous work of artificial potential functions on Euclidean spaces to the special orthogonal group.

Furthermore, we use this new configuration error function to design an adaptive update law to enable attitude convergence in the presence of uncertain disturbances. 
The stability of the proposed control systems is verified via mathematically rigorous Lyapunov analysis on $\SO$.  
In short, the proposed attitude control system in the presence of inequality constraints is geometrically exact, computationally efficient and able to handle uncertain disturbances. 
The effectiveness of this approach is illustrated via numerical simulation and demonstrated via experimental results.

\section{Problem Formulation}\label{sec:prob_form}
\subsection{Attitude Dynamics}\label{sec:att_dyn}
Consider the attitude dynamics of a rigid body.
We define an inertial reference frame and a body-fixed frame, whose origin is at the center of mass and aligned with the principle directions of the body. 
The standard orthonormal basis is denoted by \( e_i \) for \( i \in \braces{1, 2, 3} \).
The configuration manifold of the attitude dynamics is the special orthogonal group:
\begin{align*}
	\SO = \{R\in\R^{3\times 3}\,|\, R^TR=I,\;\mathrm{det}[R]=1\} ,
\end{align*}
where a rotation matrix $R\in\SO$ represents the transformation of the representation of a vector from the body-fixed frame to the inertial reference frame. 
The equations of motion are given by
\begin{gather}
	J\dot\Omega + \Omega\times J\Omega = u+W(R,\Omega)\Delta ,\label{eqn:Wdot}\\
	\dot R = R\hat\Omega ,\label{eqn:Rdot}
\end{gather}
where $J\in\R^{3\times 3}$ is the inertia matrix, and $\Omega\in\R^{3}$ is the angular velocity represented with respect to the body-fixed frame. 
The control moment is denoted by $u\in\R^{3}$, and it is expressed with respect to the body-fixed frame. 
We assume that the external disturbance is expressed by $W(R,\Omega)\Delta$, where $W(R,\Omega):\SO\times\R^{3}\rightarrow \R^{3\times p}$ is a known function of the attitude and the angular velocity.
The disturbance is represented by $\Delta\in\R^{p}$ and is an unknown, but fixed uncertain parameter.
In addition, we assume that a bound on \( W(R, \Omega) \text{ and } \Delta \) is known and given by
\begin{equation}
	\norm{W} \leq B_W , \quad \norm{\Delta} \leq B_\Delta \,, \label{eqn:W_bound}
\end{equation}
for positive constants \( B_W, B_\Delta\).
This form of uncertainty enters the system dynamics through the input channel and as a result is referred to as a matched uncertainty. 
While this form of uncertainty is easier than the unmatched variety, many physically realizable disturbances may be modeled in this manner.
For example, orbital spacecraft are subject to gravity gradient torques caused by the non-spherical distribution of mass of both the spacecraft and central gravitational body.
This form of disturbance may be represented as a body fixed torque on the vehicle.
In addition, for typical scenarios, where the spacecraft is significantly smaller than the orbital radius, the disturbance torque may be assumed constant over short time intervals.
Another terrestrial example is frequently encountered by unmanned aerial vehicles. 
Multiple actuator systems, such as quadrotor aerial vehicles, may exhibit an uneven mass distribution during load transportation or must operate in the presence of turbulence. 
Treating these effects as disturbances is a popular method to design control systems.

In~\refeqn{Rdot}, the \textit{hat} map $\wedge :\R^{3}\rightarrow\so$ represents the transformation of a vector in $\R^{3}$ to a $3\times 3$ skew-symmetric matrix such that $\hat x y = x\times y$ for any $x,y\in\R^{3}$~\cite{bullo2004}. 
More explicitly,
\begin{align*}
    \hat x = \begin{bmatrix} 0 & -x_3 & x_2 \\ x_3 & 0 & -x_1 \\ -x_2 & x_1 & 0\end{bmatrix},
\end{align*}
for $x=[x_1,x_2,x_3]^T\in\R^{3}$. 
The inverse of the hat map is denoted by the \textit{vee} map $\vee:\so\rightarrow\R^{3}$. 
We use the notation \( \hat{x}\) and \( \parenth{x}^\wedge\) interchangeably.
In particular, we use the latter form when the expression for the argument of the \textit{hat} map is complicated. 
Several properties of the hat map are summarized as
\begin{gather}
   \hat x y = x\times y = - y\times x = - \hat y x, \label{eqn:cross_product}\\
    x\cdot \hat y z = y\cdot \hat z x,\quad \hat x\hat y z = (x\cdot z) y - (x\cdot y ) z\label{eqn:STP},\\
    \widehat{x\times y} = \hat x \hat y -\hat y \hat x = yx^T-xy^T,\label{eqn:hatxy}\\
    \tr{A\hat x }=\frac{1}{2}\tr{\hat x (A-A^T)}=-x^T (A-A^T)^\vee,\label{eqn:hat1}\\
    \hat x  A+A^T\hat x=(\braces{\tr{A}I_{3\times 3}-A}x)^{\wedge}, \label{eqn:xAAx} \\
    R\hat x R^T = (Rx)^\wedge,\quad 
    R(x\times y) = Rx\times Ry\label{eqn:RxR}
\end{gather}
for any $x,y,z\in\R^{3}$, $A\in\R^{3\times 3}$ and $R\in\SO$. 
Throughout this paper, the dot product of two vectors is denoted by $x\cdot y = x^T y$ for any $x,y\in\R^n$ and the maximum eigenvalue and the minimum eigenvalue of $J$ are denoted by $\lambda_M$ and $\lambda_m$, respectively. 
The 2-norm of a matrix \( A \) is denoted by \( \norm{A} \), and its Frobenius norm is denoted by \( \norm{A} \leq \norm{A}_F = \sqrt{\tr{A^T A}} \leq \sqrt{\text{rank}(A)} \norm{A} \).

\subsection{State Inequality Constraint}

The two-sphere is the manifold of unit-vectors in \( \R^3 \), i.e., \( \Sph^2 = \{ q \in \R^3 \,  \vert \, \norm{q} = 1 \}\).
We define \( r \in \Sph^2 \) to be a unit vector from the mass center of the rigid body along a certain direction and it is represented with respect to the body-fixed frame.
For example, \( r \) may represent the pointing direction of an on-board optical sensor.
We define \( v \in \Sph^2 \) to be a unit vector from the mass center of the rigid body toward an undesired pointing direction and represented in the inertial reference frame.
For example, \( v \) may represent the inertial direction of a bright celestial object or the incoming direction of particles or other debris.
It is further assumed that optical sensor has a strict non-exposure constraint with respect to the celestial object.
We formulate this hard constraint as
\begin{align}
	r^T R^T v \leq \cos \theta , \label{eqn:constraint}
\end{align}
where we assume \( \ang{0} \leq \theta \leq \ang{90}  \) is the required minimum angular separation between \( r \) and \( R^T v \). 
The objective is to a determine a control input \( u \) that stabilizes the system from an initial attitude \( R_0 \) to a desired attitude \( R_d \) while ensuring that~\cref{eqn:constraint} is always satisfied.
\squeezeup
\section{Attitude Control on $\SO$ with Inequality Constraints}\label{sec:attitude_control}

The first step in designing a control system on a nonlinear manifold \( \Q \) is the selection of a proper configuration error function. 
This configuration error function, \( \Psi : \Q \times \Q \to \R \), is a smooth and proper positive definite function that measures the error between the current configuration and a desired configuration. 
Once an appropriate configuration error function is chosen, one can then define a configuration error vector and a velocity error vector in the tangent space \( \mathsf{T}_q \Q \) through the derivatives of \( \Psi \)~\cite{bullo2004}. 
With the configuration error function and vectors, the remaining procedure is analogous to nonlinear control design on Euclidean vector spaces. 
One chooses control inputs as functions of the state through a Lyapunov analysis on \Q.

To handle the attitude inequality constraint, we propose a new attitude configuration error function. 
More explicitly, we extend the trace form used in~\cite{bullo2004,LeeITCST13} for attitude control on \(\SO\) with the addition of a logarithmic barrier function. 
Based on the proposed configuration error function,  nonlinear geometric attitude controllers are constructed. 
A smooth control system is first developed assuming that there is  no disturbance, and then it is extended to include an adaptive update law for stabilization in the presence of unknown disturbances. 
The proposed attitude configuration error function and several properties are summarized as follows.

\begin{prop}[Attitude Error Function] \label{prop:config_error}
Define an attitude error function \( \Psi : \SO \to \R \), an attitude error vector \( e_R \in \R^3 \), and an angular velocity error vector \( e_\Omega \in \R^3 \) as follows:
\begin{gather}
	\Psi(R, R_d) = A(R, R_d) B(R) , \label{eqn:psi} \\
	e_R = e_{R_A} B(R) + A(R,R_d) e_{R_B} , \label{eqn:eR} \\
	e_\Omega = \Omega , \label{eqn:eW}
\end{gather}
with
\begin{gather}
	A(R,R_d) = \frac{1}{2} \tr{G \left( I - R_d^T R\right)} , \label{eqn:A} \\
	B(R) = 1 - \frac{1}{\alpha} \ln \left( \frac{\cos \theta -  r^T R^T v}{1 + \cos \theta}\right) . \label{eqn:B} \\
	e_{R_A} = \frac{1}{2} \parenth{G R_d^T R - R^T R_d G} ^ \vee , \label{eqn:eRA} \\
	e_{R_B} = \frac{ \left( R^T v\right)^\wedge r}{\alpha \left(r^T R^T v - \cos \theta \right)} . \label{eqn:eRB} 
\end{gather}	
where \( \alpha \in \R \) is defined as a positive constant and the matrix \( G \in \R^{3 \times 3} \) is defined as a diagonal matrix matrix for distinct, positive constants \( g_1, g_2, g_3 \in \R \).
Then, the following properties hold
\begin{enumerate}[(i)]
	\item \label{item:prop_psi_psd} \(\Psi\) is positive definite about \( R = R_d\) on $\SO$.
	\item \label{item:prop_era}The variation of \( A(R) \) with respect to a variation of \( \delta R = R \hat{\eta} \) for \( \eta \in \R^3 \) is given by
	\begin{align}\label{eq:dirDiff_A}
		\dirDiff{A}{R} &= \eta \cdot e_{R_A} ,
	\end{align}
	where the notation \( \dirDiff{A}{R} \) represents the directional derivative of $A$ with respect to $R$ along $\delta R$.
	\item \label{item:prop_erb} The variation of \( B(R) \) with respect to a variation of \( \delta R = R \hat{\eta} \) for \( \eta \in \R^3 \) is given by
	\begin{align}\label{eq:dirDiff_B}
		\dirDiff{B}{R} &= \eta \cdot e_{R_{B}}.
	\end{align}
	\item \label{item:prop_era_upbound}An upper bound of \( \norm{e_{R_A}} \) is given as:
	\begin{align}
		\norm{e_{R_A}}^2 \leq \frac{A(R)}{b_1} , \label{eqn:psi_lower_bound}
	\end{align}
	where the constant \( b_1 \) is given by \( b_1 = \frac{h_1}{h_2 + h_3} \) for 
	\begin{align*}
		h_1 &= \min\braces{g_1 + g_2, g_2 + g_3 , g_3 + g_1} ,\\
		h_2 &= \min\braces{\parenth{g_1 -g_2}^2,\parenth{g_2 -g_3}^2 , \parenth{g_3 -g_1}^2} ,\\
		h_3 &= \min\braces{\parenth{g_1 + g_2}^2, \parenth{g_2 + g_3}^2 , \parenth{g_3 + g_1}^2}.		
	\end{align*}
    \item \label{item:prop_psi_quadratic} \( \Psi \) is a locally quadratic function, which means there exist constants \( 0 < n_1 \leq n_2 \) such that
    \begin{align}\label{eq:psi_bound}
        n_1 \norm{e_R}^2 \leq \Psi(R) \leq n_2 \norm{e_R}^2 ,
    \end{align}
    on the neighborhood $D$ of the desired attitude \( R_d \)
    \begin{align}\label{eq:psi_quadratic_domain}
        D = \braces{R\in\SO  \vert \Psi < \psi < h_1, r^T R^T v < \beta < \cos \theta}
    \end{align}
    for $0<\psi < h_1 $ and $0< \beta<\cos\theta$. 
\end{enumerate}
\end{prop}
\begin{proof}
See~\Cref{proof:config_error}
\end{proof}

\Cref{eqn:psi} is composed of an attractive term, \( A (R) \) toward the desired attitude, and a repulsive term, \( B(R) \) away from the undesired direction \( R^T v \).
In order to visualize the attitude error function on \( \SO \), we utilize a spherical coordinate representation.
Recall, that the spherical coordinate system represents the position of a point relative to an origin in terms of a radial distance, azimuth, and elevation.
This coordinate system is commonly used to define locations on the Earth in terms of a latitude and longitude.
Similarly, the positions of celestial objects are defined on the celestial sphere in terms of right ascension and declination. 
\begin{figure}[htbp]%
    \centering 
    \subcaptionbox{ Attractive \( A(R) \) \label{fig:attract_error} }{\includegraphics[trim={8mm 2mm 8mm 2mm}, clip, width=0.3\textwidth]{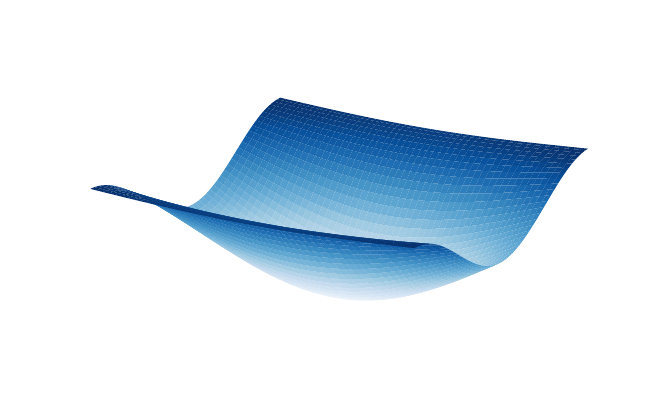}}~
    \subcaptionbox{ Repulsive \( B(R) \) \label{fig:avoid_error} }{\includegraphics[trim={8mm 2mm 8mm 2mm}, clip, width=0.3\textwidth]{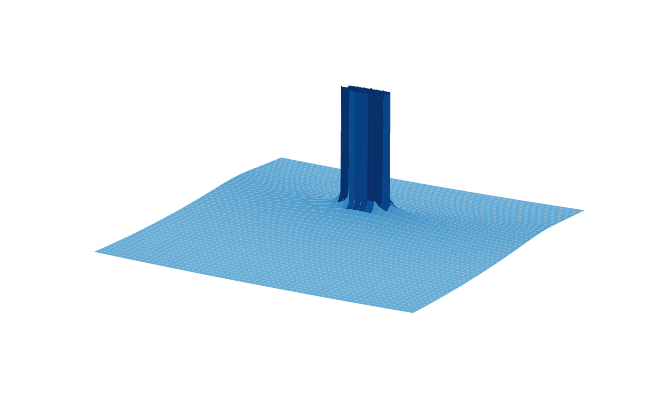}}~
    \subcaptionbox{ Combined \( \Psi \) \label{fig:combined_error} }{\includegraphics[trim={8mm 2mm 8mm 2mm}, clip, width=0.3\textwidth]{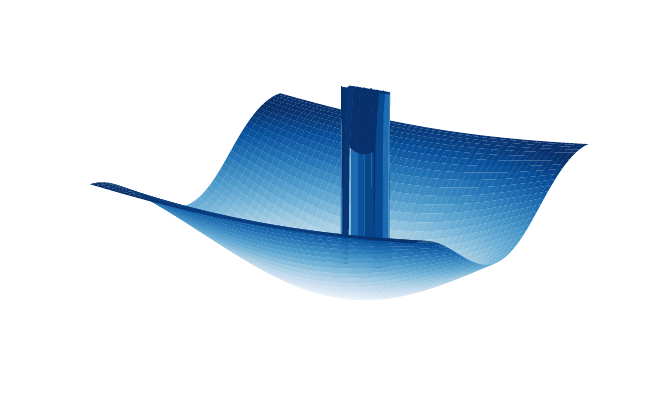}}%
    \caption{Visualization of Configuration Error Functions using spherical coordinate representation}
    \label{fig:config_error} 
\end{figure}%
We apply this concept and parameterize the rotation matrix \( R \in \SO \) in terms of the spherical angles \( \SI{-180}{\degree} \leq \lambda \leq \SI{180}{\degree}  \) and \( \SI{-90}{\degree} \leq \beta \leq \SI{90}{\degree} \). 
Using the elementary Euler rotations, the rotation matrix is now defined as \( R = \exp( \lambda \hat{e}_2) \exp( \beta \hat{e}_3) \).
We iterate over the domains of \( \lambda\) and \(\beta\) in order to rotate the body-fixed vector \( r \) throughout the two-sphere \( \S^2 \).
Applying this method,~\cref{fig:config_error} allows us to visualize the error function on \( \SO \).
The horizontal axes of~\cref{fig:config_error} represent the domain of the spherical angles \( \lambda \) and \( \beta \) in degrees, while the vertical axes represent the unitless magnitude of the error functions defined in~\cref{eqn:psi,eqn:A,eqn:B}.
The attractive error function, given by~\cref{eqn:A}, has been previously used for attitude control on \(\SO\).
The potential well of \( A(R)\) is illustrated in~\cref{fig:attract_error}, where the desired attitude lies at the minimum of \( A(R) \).

To incorporate the state inequality constraints we apply a logarithmic barrier term.
Barrier functions are typically used in optimal control and motion planning applications.
A visualization of the repulsive error function is presented in~\cref{fig:avoid_error} which shows that as the boundary of the constraint is neared, or \( r^T R^T v \to \cos \theta \), the barrier term increases, \( B \to \infty\).
We use the scale factor~\(\frac{1}{1+\cos \theta} \) to ensure that \( \Psi \) remains positive definite.
The logarithmic function is popular as it quickly decays away from the constraint boundary.
The positive constant \( \alpha \) serves to shape the barrier function.
As \( \alpha \) is increased the impact of \( B(R) \) is reduced away from the constraint boundary. 
The superposition of the attractive and repulsive functions is shown in~\cref{fig:combined_error}.
The control system is defined such that the attitude trajectory follows the negative gradient of \( \Psi \) toward the minimum at \( R = R_d \), while avoiding the constrained region.

While~\cref{eqn:B} represents a single inequality constraint given as~\cref{eqn:constraint}, it is readily generalized to multiple constraints of an arbitrary orientation. 
For example, the configuration error function can be formulated as $\Psi=A[1+\sum_i C_i]$, where $C_i$ has the form of $C_i=B-1$ for the $i$-th constraint. 
In this manner, one may enforce multiple state inequality constraints, and we later demonstrate this through numerical simulation. 
This is in contrast to many previous approaches which are computationally difficult to extend to situations with multiple constraints.
We present the dynamics of the configuration error function in~\Cref{prop:error_dyn}, which are used in the subsequent development of the nonlinear control system.

\begin{prop}[Error Dynamics]\label{prop:error_dyn}
	The attitude error dynamics for \( \Psi, e_R, e_\Omega \) satisfy 
	\begin{gather}
		\diff{}{t} \parenth{\Psi} = e_R \cdot e_\Omega , \label{eqn:psi_dot}\\
		\diff{}{t} \parenth{e_R} = \dot{e}_{R_A} B + e_{R_A} \dot{B} + \dot{A}e_{R_B} + A \dot{e}_{R_B} , \label{eqn:eR_dot} \\
		\diff{}{t} \parenth{e_{R_A}} = E(R, R_d) e_\Omega , \label{eqn:eRA_dot} \\
		\diff{}{t} \parenth{e_{R_B}} = F(R) e_\Omega , \label{eqn:eRB_dot} \\
    	\diff{}{t} \parenth{A(R)} = e_{R_A} \cdot e_\Omega , \label{eqn:A_dot} \\
		\diff{}{t} \parenth{B(R)} = e_{R_B} \cdot e_\Omega , \label{eqn:B_dot} \\
		\diff{}{t} \parenth{e_\Omega} = J^{-1} \parenth{-\Omega \times J \Omega + u + W(R, \Omega) \Delta} , \label{eqn:eW_dot}
	\end{gather}
	where the matrices \(E(R,R_d), F(R) \in \R^{3\times3} \) are given by
	\begin{align}
		E(R,R_d) = &\frac{1}{2} \parenth{\tr{R^T R_d G}I - R^T R_d G} , \label{eqn:E} \\
		F(R) = &\frac{1}{\alpha \parenth{r^T R^T v - \cos \theta}} \left[\parenth{v^T R r} I - R^T v r^T + \right. \nonumber \\
		& \left. \frac{R^T \hat{v} R r v^T R \hat{r}}{\parenth{r^T R^T v - \cos \theta}}\right] . \label{eqn:F}
	\end{align}
\end{prop}
\begin{proof}
See~\Cref{proof:error_dyn}
\end{proof}

\begin{figure}[t]
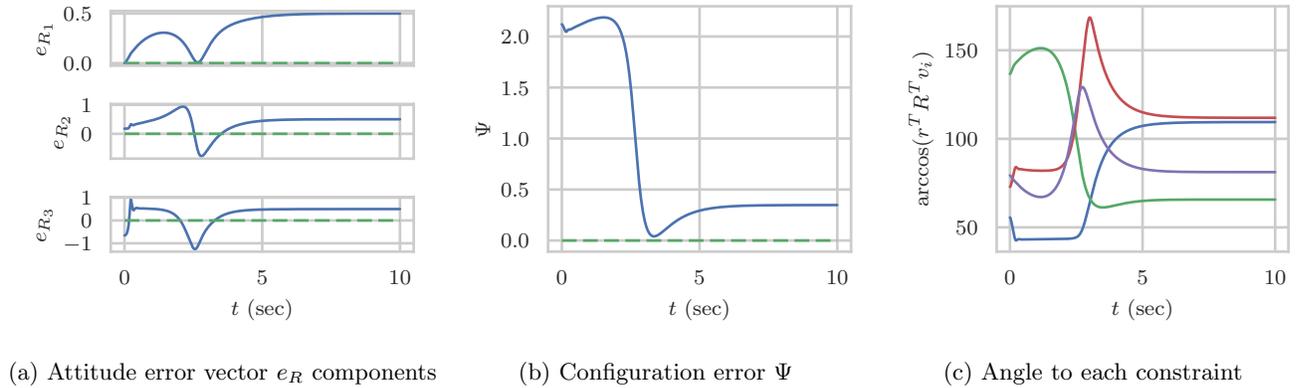

    \centering 
    \subcaptionbox{Attitude error vector \(e_R\) components \label{fig:eR_con}}{\input{pgf/eR_noadapt.pgf}}~
    \subcaptionbox{Configuration error \( \Psi \) \label{fig:Psi_con}}{\input{pgf/Psi_noadapt.pgf}}~
   \subcaptionbox{Angle to each constraint \label{fig:con_angles_con}}{\input{pgf/ang_con_noadapt.pgf}}~
    \caption{Attitude stabilization without adaptive update law}
    \label{fig:con} 
\end{figure}
\subsection{Attitude Control without Disturbance}
We introduce a nonlinear geometric controller for the attitude stabilization of a rigid body.
We first assume that there is no disturbance, i.e., \( \Delta = 0 \), and present a nonlinear controller in~\Cref{prop:att_control}.
\begin{prop}[Attitude Control]\label{prop:att_control}
	Given a desired attitude command \( \parenth{R_d, \Omega_d = 0} \), which satisfies the constraint~\cref{eqn:constraint}, and positive constants \( k_R, k_\Omega \in \R \) we define a control input \( u \in \R^3 \) as follows
	\begin{gather}
		u = -k_R e_R - k_\Omega e_\Omega + \Omega \times J \Omega . \label{eqn:nodist_control}
	\end{gather}
	Then the zero equilibrium of the attitude error is asymptotically stable, and the inequality constraint is always satisfied.
\end{prop}
\begin{proof}
See~\Cref{proof:att_control}
\end{proof}

\subsection{Adaptive Control}
We extend the results of the previous section with the addition of a fixed but unknown disturbance \( \Delta \).
This scenario is typical of many mechanical systems and represents unmodelled dynamics or external moments acting on the system.
For example, Earth orbiting spacecraft typically experience a torque due to a gravitational gradient.
Aerial vehicles will similarly experience external torques due to air currents or turbulence.
An adaptive control system is introduced to asymptotically stabilize the system to a desired attitude while ensuring that state constraints are satisfied. 

\begin{prop}[Bound on \( \bm{\dot{e}_R} \)]\label{prop:eR_dot_bound}
Consider the neighborhood \( D \), given in~\cref{prop:config_error}, about the desired attitude, then the following statements hold:
\begin{enumerate}[(i)]
    \item \label{item:prop_eR_dot_bound_AB} Upper bounds of \( A(R) \) and \( B(R) \) are given by
        \begin{gather}
            \norm{A} < b_2 \norm{e_{R_A}}^2 < c_A  , \quad \norm{B} < c_B . \label{eqn:AB_bound}
        \end{gather}
        where the constant \( b_2\) is given by \( b_2 = \frac{h_1 h_4}{h_5 \parenth{h_1 - \psi}}\) for
        \begin{align*}
            h_4 &= \min\braces{g_1 + g_2, g_2 + g_3 , g_3 + g_1} ,\\
            h_5 &= \min\braces{\parenth{g_1 + g_2}^2,\parenth{g_2 + g_3}^2 , \parenth{g_3 + g_1}^2} .\\
        \end{align*}
    \item \label{item:prop_eR_dot_bound_EF} Upper bounds of \( E(R,R_d) \) and \( F(R) \) are given by
        \begin{gather}
            \norm{E} \leq \frac{1}{\sqrt{2}} \tr{G}  , \label{eqn:E_bound} \\
            \norm{F} \leq \frac{\parenth{\beta^2 + 1}\parenth{\beta - \cos \theta}^2 + 1 + \beta^2 \parenth{\beta^2-2}}{\alpha^2 \parenth{\beta-\cos \theta}^4} . \label{eqn:F_bound}
        \end{gather}
    \item Upper bounds of the attitude error vectors \( e_{R_A} \) and \( e_{R_B} \) are given by
        \begin{gather}
            \norm{e_{R_A}} \leq \sqrt{\frac{\psi}{b_1}}, \label{eqn:eRA_bound} \\
            \norm{e_{R_B}} \leq \frac{\sin\theta}{\alpha \parenth{\cos \theta - \beta}}. \label{eqn:eRB_bound}
        \end{gather}
\end{enumerate}
These results are combined to yield a maximum upper bound of the time derivative of the attitude error vector \( \dot{e}_R \) as
\begin{gather}
	\norm{\dot{e}_R} \leq H \norm{e_\Omega} ,\label{eqn:eR_bound}
\end{gather}
where  \( H \in \R \) is defined as
\begin{gather}
	H = \norm{B} \norm{E} + 2 \norm{e_{R_A}} \norm{e_{R_B}} + \norm{A}\norm{F}. \label{eqn:H}
\end{gather}
\end{prop}
\begin{proof}
See~\Cref{proof:eR_dot_bound}
\end{proof}

Adaptive control is typically used in dynamical systems with varying or uncertain components.
In~\Cref{prop:adaptive_control}, we present an adaptive attitude controller which handles uncertain disturbances while satisfying the state inequality constraints.

\begin{prop}[Adaptive Attitude Control]\label{prop:adaptive_control}
Given  a desired attitude command \( (R_d, \Omega_d = 0 )\) and positive constants \( k_R, k_\Omega, k_\Delta, c \in \R \), we define a control input \( u \in \R^3\) and an adaptive update law for the estimated uncertainty \( \bar{\Delta} \) as follows:
\begin{align}
	u &= - k_R e_R - k_\Omega e_\Omega + \Omega \times J \Omega - W \bar{\Delta} , \label{eqn:adaptive_control} \\
	\dot{\bar{\Delta}} &= k_\Delta W^T \parenth{e_\Omega + c e_R} . \label{eqn:delta_dot}
\end{align}
If \( c \) is chosen such that
\begin{gather}
	0 < c < \min \braces{\sqrt{\frac{2 \lambda_m k_R n_1}{\lambda_M^2}},
	\frac{4 k_R k_\Omega}{k_\Omega^2 + 4 k_R \lambda_M H}} , \label{eqn:c_bound}
\end{gather}
  the zero equilibrium of the error vectors is stable in the sense of Lyapunov. Furthermore, $e_R,e_\Omega\rightarrow 0$ as $t\rightarrow\infty$, and $\bar\Delta$ is  bounded.
\end{prop}
\begin{proof}
See~\Cref{proof:adaptive_control}
\end{proof}

Nonlinear adaptive controllers have been developed for attitude stabilization in terms of modified Rodriguez parameters and quaternions, as well as attitude tracking in terms of Euler angles. 
The proposed control system is developed on \(\SO\) and avoids the singularities of Euler angles and Rodriguez parameters while incorporating state inequality constraints. 
In addition, the unwinding and double coverage ambiguity of quaternions are also completely avoided. 
The control system handles uncertain disturbances while avoiding constrained regions.

Compared to the previous work on constrained attitude control, we present a geometrically exact control system without parameterizations.
The controller is designed on the true configuration manifold, \( \SO \), and is free from the issues associated with other attitude representations.
In addition, we incorporate state inequality constraints for the first time on \( \SO \).
The presented control system is computed in real-time and offers significant computational advantages over previous iterative methods. 
In addition, the rigorous mathematical proof guarantees stability.
This is in contrast to many of the previous methods which offer no stability guarantees.
The presented analysis offers provable bounds on the expected motion, which are critical for hardware implementation or mission critical applications.

\section{Numerical Examples}\label{sec:numerical_simulation}

\begin{figure} 
  \centering 
  \subcaptionbox{ Configuration error \( \Psi \)\label{fig:Psi_adapt} }{\input{pgf/Psi_adapt.pgf}}~
  \subcaptionbox{ Angle to each constraint \label{fig:con_angles} }{\input{pgf/ang_con_adapt.pgf}}\\
  \subcaptionbox{ Disturbance estimate \( \bar \Delta \) components \label{fig:delta_adapt} }{\input{pgf/dist_adapt.pgf}}~
    \subcaptionbox{ Attitude trajectory \label{fig:cad_adapt} }{\includegraphics[trim={10cm 0 10cm 0},clip,width=0.4\columnwidth]{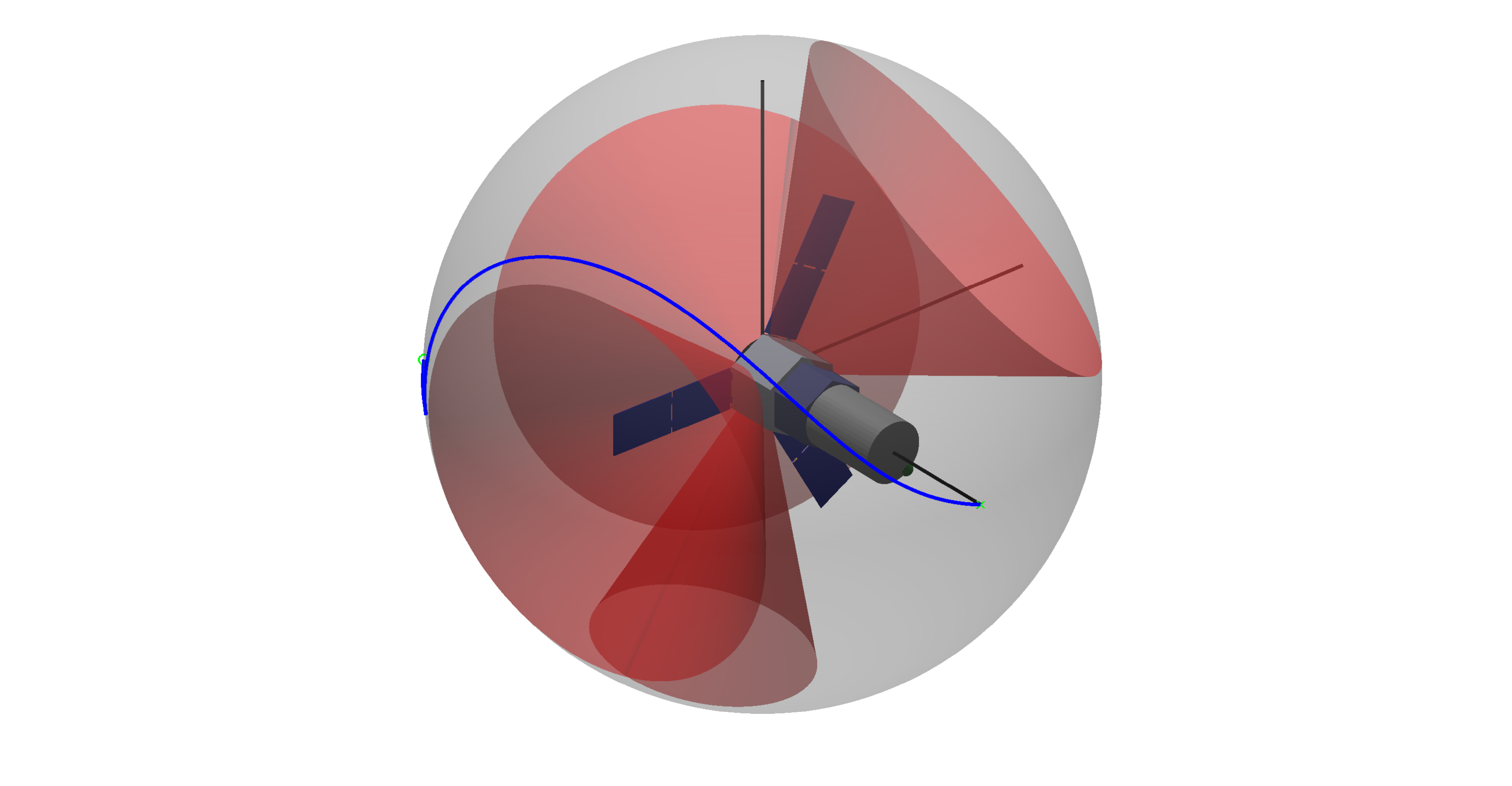} }
  \caption{Attitude stabilization with adaptive update law}
  \label{fig:adapt} 
\end{figure}

We demonstrate the performance of the proposed control system via numerical simulation.
The inertia tensor of a rigid body is given as
\begin{gather*}
    J = \begin{bmatrix}
	\num{5.5} & \num{0.06} & \num{-0.03} \\
	\num{0.06} & \num{5.5} & \num{0.01} \\
	\num{-0.03} & \num{0.01} & \num{0.1}
    \end{bmatrix} \times 10^{-3}~\si{\kilo\gram\meter\squared} .
\end{gather*} 
The control system parameters are chosen as
\begin{gather*}
	G = \text{diag} [0.9,1.1,1.0], \quad k_R = 0.4 , \quad	k_\Omega = 0.296 ,\\
	c = 1.0 , \quad k_\Delta = 0.5 , \quad \alpha = 15 .
\end{gather*}
The diagonal matrix \( G \) serves as a weighting matrix for the relative difference between \( R \) and \( R_d \). 
Using this term, the control designer can modify the shape of the attractive error function, given in~\cref{eqn:A}, and the resulting behavior of the closed loop system.
Similarly, the constant \( \alpha \) is used to modify the shape of the repulsive error function, given in~\cref{eqn:B}.
In general, this term is derived from the system design and the nature of the obstacles in the dynamic environment.
For example, a system wishing to avoid pointing at a diffuse obstacle, such as incoming debris, may chose an appropriate value of \( \theta \), based on the best available information, and a relatively low \( \alpha \) to ensure additional safety margin near the constraint boundary. 
Similarly, in an environment with several densely spaced obstacles, such as that presented in~\cref{fig:cad_adapt}, a much larger \( \alpha \) would enable more aggressive maneuvers which pass closer to the constraint boundary without violation.
This would increase the allowable region of operation in a highly constrained environment.
The parameters \( k_R, k_\Omega, c, k_\Delta\) are control parameters used to modify the closed-loop behavior of the system.
It is straightforward to chose \( k_R, k_\Omega, k_\Delta\), using a linear analysis, to satisfy desired response criteria, such as settling time or percent overshoot~\cite{nise2004}.
A body fixed sensor is defined as \(r = [1,0,0]\), while multiple inequality constraints are defined in~\Cref{tab:constraints}.
The simulation parameters are chosen to be similar to those found in~\cite{lee2011b}, however we increase the size of the constraint regions to create a more challenging scenario for the control system.
\vspace{-2.5mm}
\begin{table}[htbp]
\caption{Constraint Parameters~\label{tab:constraints}}
\begin{center}\begin{tabular}{lc}
Constraint Vector (\( v \)) & Angle (\( \theta \)) \\ \hline \hline 
\([0.174,\,-0.934,\, -0.034]^T\) & \ang{40} \\ \hline 
\([0 ,\, 0.7071 ,\, 0.7071]^T\) & \ang{40} \\ \hline 
\([-0.853 ,\, 0.436 ,\, -0.286]^T\) & \ang{40} \\ \hline 
\([-0.122 ,\,-0.140,\, -0.983]^T\) & \ang{20}\end{tabular} 
\end{center}
\end{table}
\vspace{-2.5mm}
The initial state is defined as \(R_0 =  \exp(\ang{225} \times \frac{\pi}{180} \hat{e}_3), \Omega_0 = 0\), with \( e_3 = \begin{bmatrix} 0 & 0 & 1 \end{bmatrix}^T \).
The desired state is \( R_d = I,\Omega_d = 0\).
We show simulation results for the system stabilizing about the desired attitude with and without the adaptive update law from~\Cref{prop:adaptive_control}.
We assume a fixed disturbance of \(\Delta = \begin{bmatrix} 0.2 & 0.2 & 0.2 \end{bmatrix}^T \si{\newton\meter}\), with the function \( W(R,\Omega) = I \).
This form is equivalent to an integral control term which penalizes deviations from the desired configuration.
The first term of~\cref{eqn:delta_dot} has the effect of increasing the proportional gain of the control system, since the time derivative of the attitude error vector, \( \dot{e}_{R} \), is linear with respect to the angular velocity error vector \( e_\Omega\).

Simulation results without the adaptive update law are shown in~\cref{fig:con}.
\Cref{fig:eR_con} shows each component of the attitude error vector,~\cref{eqn:eR}, over the simulation time span.
\Cref{fig:Psi_con} shows the  magnitude of the combined error function,~\cref{eqn:psi}.
Without the update law, the system does not achieve zero steady state error. 
\Cref{fig:Psi_con} shows that the configuration error function does not converge to zero and there exist steady state errors.
In spite of the uncompensated disturbance, the system is able to avoid the constrained regions as shown in~\cref{fig:con_angles_con}.
The angle to each of the constraints, which is measured in degrees and given by \( \arccos(r^T R^T v_i) \), is always greater than the specified angle, \( \theta_i \), in~\Cref{tab:constraints}.

\Cref{fig:adapt} shows the results with the addition of the adaptive update law.
\Cref{fig:Psi_adapt,fig:con_angles} are equivalent to~\cref{fig:Psi_con,fig:con_angles_con} with the exception of the addition of the adaptive update law.
The addition of the adaptive update law allows the system to converge to the desired attitude in the presence of constraints.
The path of the body fixed sensor in the inertial frame, namely \( R r \), is illustrated in~\cref{fig:cad_adapt} by the blue trajectory.
The rendering of the spacecraft is presented in the desired, or final, orientation of the simulation.
The inequality constraints from~\Cref{tab:constraints} are depicted as red cones, where the cone half angle is \( \theta \).
The control system is able to asymptotically converge to the desired attitude.
\Cref{fig:con_angles} shows that the angle, \( \arccos(r^T R^T v_i) \) and measured in degrees, between the body fixed sensor and each constraint is satisfied for the entire maneuver.
In addition, the estimate of the disturbance converges to the true value as shown in~\cref{fig:delta_adapt}.

Both control system are able to automatically avoid the constrained regions. 
In addition, these results show that it is straightforward to incorporate an arbitrary amount of large constraints.
In spite of this challenging configuration space, the proposed control system offers a simple method of avoiding constrained regions.
These closed-loop feedback results are computed in real time and offer a significant advantage over typical open-loop planning methods.
These results show that the proposed geometric adaptive approach is critical to attitude stabilization in the presence of state constraints and disturbances.

\subsection{Attitude Parameterizations}\label{ssec:attitude_parameterization}
Attitude parameterizations, such as Euler angles and Quaternions, are frequently used in the aerospace and astrodynamics communities~\cite{vallado2007}.
For example, Euler angle sequences are frequently used to describe the transformation between a variety of reference frames used to describe the position and orientation of the orbit of Earth satellites~\cite{vallado2007}.
In addition, quaternions were used during the operation of Skylab and the NASA Space Shuttle~\cite{hughes2004}.
However, the choice of attitude parameterization plays a critical role in control design and the resulting motion of the system.

Euler angle sequences are a minimum, three-parameter set of angles which describe the transformation between two reference frames.
Using Euler angles, we can represent any general rotation as a sequence of three intermediate rotations~\cite{shuster1993}.
By convention, there are \num{24} possible Euler angle sequences for any given rotation.
In addition, Euler angles are a minimum representation, as only three angles, and the associated sequence, are required to describe the three angular degrees of freedom of the rigid body.
However, there is great ambiguity in the representation of the attitude as there are many equivalent Euler angle sequences for a given attitude of the system.
Therefore, great care must be taken in the control system design to ensure that a consistent sequence is used. 
Furthermore, it has been shown that no minimal attitude representation can describe orientations both globally and without singularities~\cite{hughes2004,bhat2000}.
These singularities can cause significant difficulties during control design and hardware implementation.

To demonstrate the effect of the kinematic singularities inherent with Euler angles we will represent the attitude of the body fixed reference frame, \( \vecbf{b}_i \), with respect to the inertial frame, \( \vecbf{e}_i\), in terms of the 3-1-3 Euler angle sequence.
More explicitly, this corresponds to the rotation sequence \( \theta_1 \vecbf{b}_3 , \theta_2 \vecbf{b}_1, \theta_3 \vecbf{b}_3 \).
The rotation matrix, \( R(\theta_1, \theta_2, \theta_3) \), corresponding to this sequence is 
\begin{align}\label{eq:euler313}
    \begin{bmatrix}
        -s_1 c_2 s_3 + c_3 c_1 & -s_1 c_2 c_3 - s_3 c_1 & s_1s_2 \\
        c_1 c_2 s_3 + c_3 s_1 & c_1 c_2 c_3 - s_3 s_1 & - c_1 s_2 \\
        s_2 s_3 & s_2 c_3 & c_2
    \end{bmatrix} ,
\end{align}
where \( s_i, c_i \) represent \( \sin \theta_i, \cos \theta_i \) for \( i = \braces{1,2,3}\).
Using this representation, the kinematic differential equations for the associated Euler angles are given as
\begin{align}\label{eq:euler313_diff}
    \begin{bmatrix}
        \dot{\theta}_1 \\ \dot{\theta}_2 \\ \dot{\theta}_3 
    \end{bmatrix}
    =
    \begin{bmatrix}
        \slfrac{\parenth{\Omega_1 s_3 + \Omega_2 c_3}}{s_2} \\
        \Omega_1 c_3 - \Omega_2 s_3 \\
        -\slfrac{\parenth{\Omega_1 s_3 + \omega_2 c_3}c_2}{s_2} + \Omega_3
    \end{bmatrix} .
\end{align}
From~\cref{eq:euler313_diff}, it is immediately clear that a singularity exists when \( \sin \theta_2 = 0 \) or equivalently, \( \theta_2 = 0, \pm \pi \). 
In the vicinity of the singularity, the angular velocities of the Euler angles will tend to approach \( \pm \infty \) and the angular velocities will experience instantaneous sign changes.
Furthermore, all Euler angle sequences will exhibit a similar singularity at either \( \theta_2 = 0, \pm \pi \) or \( \theta_2 = \pm \frac{\pi}{2}, \pm \frac{3\pi}{2} \).
Therefore simply switching the sequence does not alleviate the issue, but rather only moves the singularity.
As a result, Euler angles are not appropriate for systems which experience large angular rotations, such as those demonstrated in~\cref{fig:adapt}, or control systems which rely on the angular velocities \( \theta_i \).
\subsection{Time-varying Disturbance}\label{ssec:time_varying}
The form of the uncertainty, given in~\cref{eqn:Wdot}, is commonly used in the adaptive control literature~\cite{LeeITCST13,ioannou2012}. 
A wide variety of realistic disturbances, such as gravitational gradients or malfunctioning thrusters for spacecraft scenarios, are accurately represented via this model. 
In addition, it is possible to represent the uncertainty of a time-varying inertia matrix as an equivalent external disturbance. 
For example, Euler's law gives the relationship for the rate of change of angular momentum as
\begin{align*}
    M_{ext} = \dot{\vecbf{H}} = \dot{J} \vecbf{\Omega} + J \dot{\vecbf{\Omega}} .
\end{align*}
Using this, we can see that an instantaneous change in \( J \) is proportional to an external moment.
Finally, it has been shown that this adaptive control formulation is able to handle time-varying disturbances under some mild assumptions~\cite{ioannou2012}. 

We demonstrate the ability to handle an uncertain time-varying disturbance via numerical example.
The system is identical to the one presented in~\Cref{sec:numerical_simulation}, however we modify the external disturbance. 
The external disturbance is the superposition of constant and time-varying terms as
\begin{align*}
    \Delta = \begin{bmatrix} 0.2 \\ 0.2 \\0.2 \end{bmatrix} + 0.02 \begin{bmatrix} \sin 9 t \\ \cos 9 t \\ \frac{1}{2} \parenth{\sin 9t + \cos 9t}\end{bmatrix} \si{\newton\meter}.
\end{align*}
We define a constraint in the inertial frame as \( v = [\frac{1}{\sqrt{2}}, \frac{1}{\sqrt{2}}, 0]^T \) with \( \theta = \ang{12} \).
The initial state is defined as \(R(0) = \exp( \frac{\pi}{2} \hat{e}_3) \), while the desired state is \(R_d =I \).
The goal is to rotate the vehicle about the \( e_3 \) axis while avoiding the obstacle and compensating for the time-varying disturbance. 

\begin{figure}[htbp]
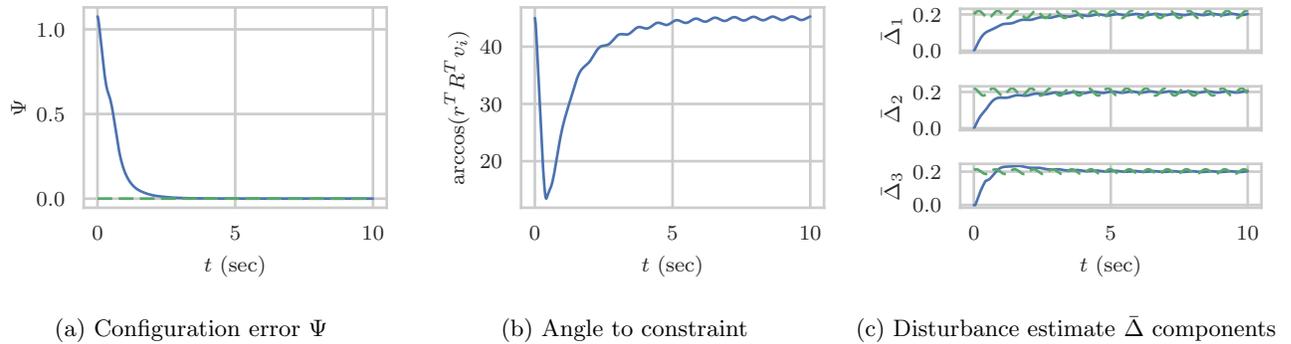

    \centering 
    \subcaptionbox{ {Configuration error \( \Psi \)} \label{fig:Psi_tv} }{\input{pgf/Psi_timevarying.pgf}}~
    \subcaptionbox{ {Angle to constraint } \label{fig:con_angle_tv} }{\input{pgf/ang_con_timevarying.pgf}}~
    \subcaptionbox{ {Disturbance estimate \(\bar\Delta\) components} \label{fig:Delta_tv} }{\input{pgf/dist_timevarying.pgf}}
    \caption{Time-varying external disturbance simulation}
    \label{fig:tv} 
\end{figure}

\Cref{fig:tv} demonstrates the ability for the adaptive controller, which is presented in~\Cref{prop:adaptive_control}, to handle time-varying disturbances.
\Cref{fig:Psi_tv} shows the non-dimensional value of the configuration error function and demonstrates that the adaptive controller is able to stabilize the system to the desired attitude configuration.
In addition,~\cref{fig:con_angle_tv} shows that the constraint is never violated as the angle between the body-fixed sensor \( r \) and the constraint \( v \) is greater than \SI{12}{\degree} over the entire attitude maneuver.
We can see in~\cref{fig:Delta_tv} that that estimate \(\bar \Delta \) for each of the components accurately tracks the true disturbance after approximately~\SI{5}{\second}.

\section{Experiment on Hexrotor UAV}\label{sec:experiment}

\begin{figure}[htbp]
    \centering
    \includegraphics[width = 0.50\columnwidth]{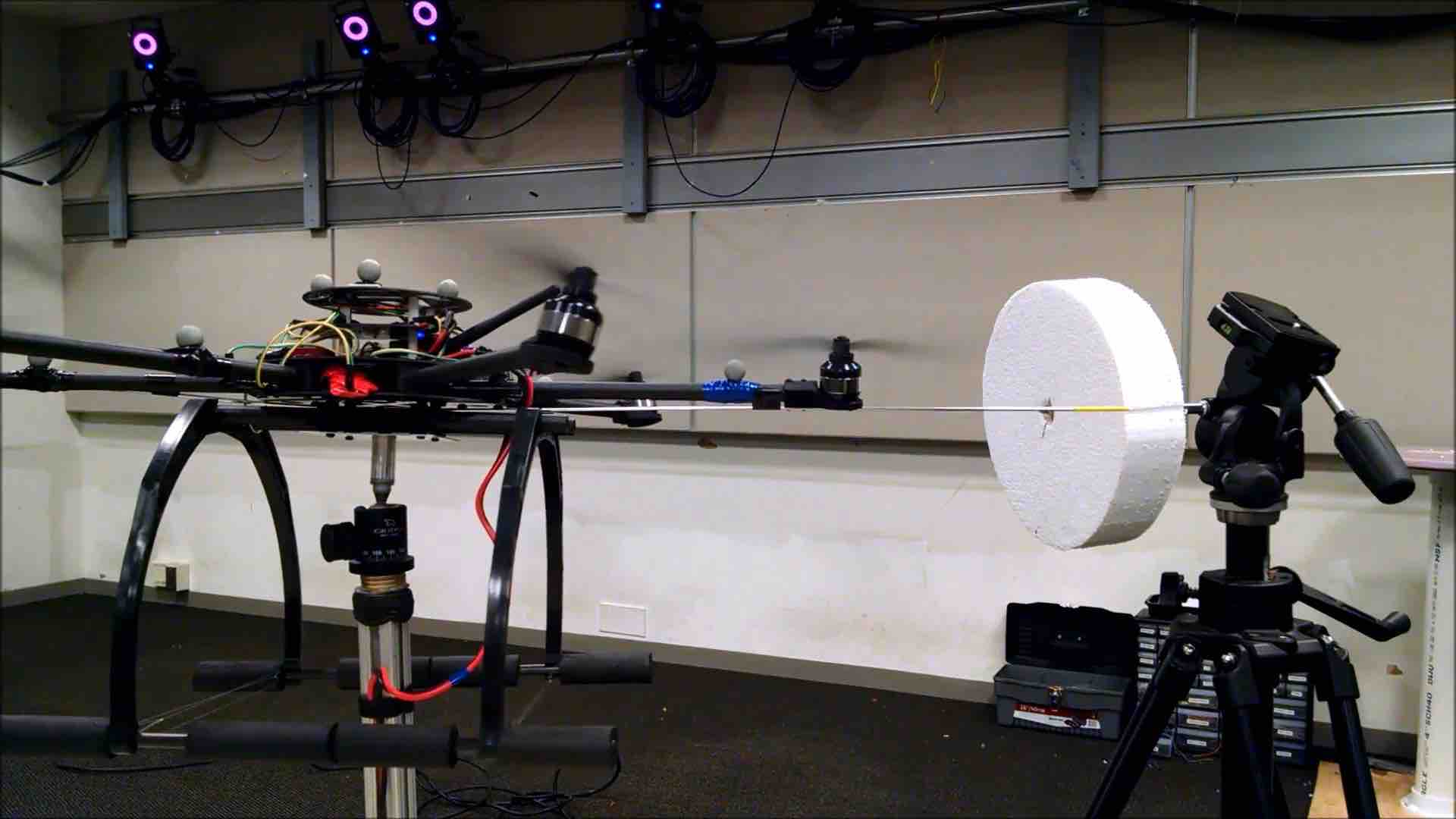}
    \caption{Attitude control testbed~\label{fig:hexrotor}}
\end{figure}
\vspace{-3mm}
A hexrotor unmanned aerial vehicle (UAV), as seen in~\cref{fig:hexrotor}, has been developed at the Flight Dynamics and Controls Laboratory (FDCL) at the George Washington University~\cite{kaufman2014}.
The UAV is composed of three pairs of counter-rotating propellers. 
Typical UAVs are composed of four or more co-planar propellers.
As a result, these systems are underactuated and unable to impart a force along every degree of freedom.
For example, quadrotor UAVs are unable to translate laterally without first conducting a rotation.
Conversely, the propeller pairs of the hexrotor are angled relative to one another to allow for a fully actuated rigid body.
This allows the hexrotor to impart a force in any direction and a moment about any axis. 
\begin{figure}[t]
    \centering 
    \subcaptionbox{ {Attitude error vector \(e_R\)} components\label{fig:eR_exp} }{\input{pgf/eR_exp.pgf}}~
    \subcaptionbox{ {Configuration error \( \Psi \)} \label{fig:Psi_exp} }{\input{pgf/Psi_exp.pgf}}~
    \subcaptionbox{ {Control input \( u\) } components\label{fig:u_exp} }{\input{pgf/u_exp.pgf}}\\
    \subcaptionbox{ {Attitude Trajectory } \label{fig:traj_exp} }{\includegraphics[width=0.3\columnwidth]{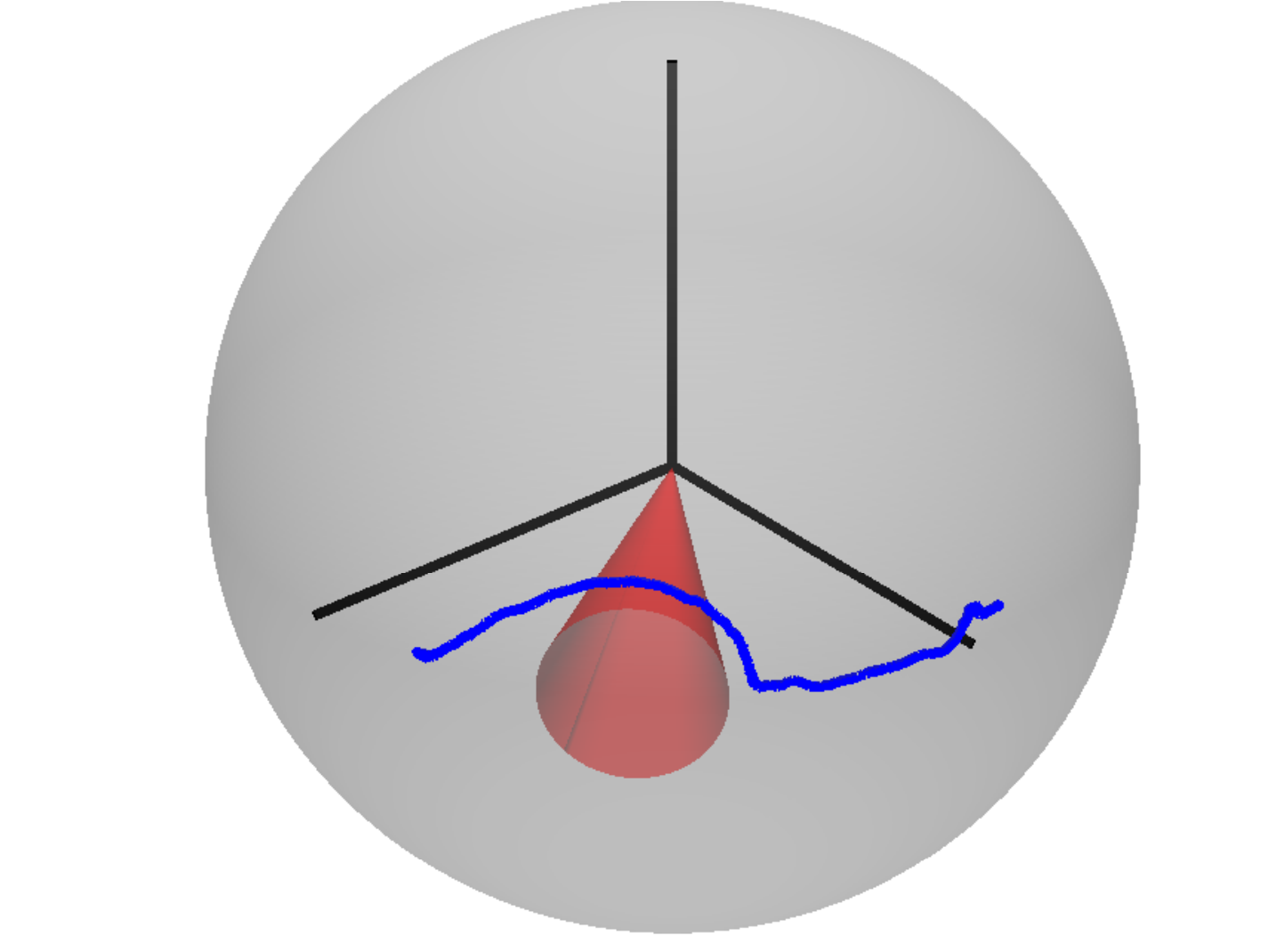} }~
    \subcaptionbox{ {Angle to constraint } \label{fig:con_angle_exp} }{\input{pgf/ang_con_exp.pgf}}
    \caption{Constrained Attitude stabilization experiment}
    \label{fig:exp} 
\end{figure}
Attitude information is measured by a combination of both on and off board sensor systems.
The VectorNav VN-100 is a rugged, miniature high-performance inertial measurement unit which provides high frequency angular velocity measurements.
A Vicon motion capture system is installed within the test environment and used to provide high accuracy attitude measurements. 
A series of reflective markers are placed on the hexrotor and their relative position is captured by a series of infrared optical cameras. 
Assuming a fixed rigid body, the Vicon system is able to derive the attitude of the hexrotor and transmit this data to the processor onboard the hexrotor.
The control input is computed on-board, using the full state measurement, and implemented at approximately \SI{100}{\hertz}.
In order to constrain the motion, allowing us to test only the attitude dynamics, we attach the hexrotor to a spherical joint.
Since, the center of rotation is below the center of gravity of the hexrotor there is a destabilizing gravitational moment.
The resulting attitude dynamics are similar to an inverted pendulum model.
We augment the control input in~\cref{eqn:adaptive_control} with an additional term to negate the effect of the gravitational moment.

A sensor pointing direction is defined in the body-fixed frame of the hexrotor as \( r = [1,0,0]^T \).
We define an obstacle in the inertial frame as \( v = [\frac{1}{\sqrt{2}}, \frac{1}{\sqrt{2}}, 0]^T \) with \( \theta = \ang{12} \).
An initial state is defined as \(R(0) = \exp( \frac{\pi}{2} \hat{e}_3) \), while the desired state is \(R_d =I \).
This results in the UAV performing a \ang{90} yaw rotation about the vertical axis of the spherical joint and the constrained region is on the shortest path connecting $R_0$ and $R_d$. 
The attitude control system is identical to the one presented in~\Cref{prop:adaptive_control} with the exception of a gravity moment term, \( M_g = r_{cg} \times m g R^T e_3\) which represents the gravitational moment due to the distance \( r_{cg} \) between the center of mass and the center of rotation. 
In addition, the following parameters were also modified: \(k_R = 0.4, k_\Omega = 0.7 ,c = 0.1 , \alpha = 8 \text{ and } k_\Delta = 0.05\) to account for the differences in the hardware model of the hexrotor.

The experimental results are shown in~\Cref{fig:exp}.
\Cref{fig:eR_exp} shows the behavior of each of the components of the attitude error vector, defined by~\cref{eqn:eR}, over the experiment time span.
\Cref{fig:Psi_exp} shows the time history of the attitude error function, defined by~\cref{eqn:psi}.
\Cref{fig:u_exp} shows the magnitude of each component of the control input in \si{\newton\meter}, which is computed from~\cref{eqn:adaptive_control}.
Finally,~\cref{fig:con_angle_exp} shows the angle between the body-fixed sensor and the obstacle in degrees.
In order to maneuver the system ``close" to the constrained zone we utilize several intermediary set points on either side of the obstacle.
From the initial attitude the hexrotor rotates to the first set point, pauses, and then continues around the obstacle to the second set point before continuing toward the desired attitude.
As a result this creates the stepped behavior of the configuration error history as shown in~\cref{fig:Psi_exp}.

The hexrotor avoids the constrained region illustrated by the circular cone in~\cref{fig:traj_exp}, by rotating around the boundary of the constraint. 
\Cref{fig:con_angle_exp} shows the angle, \( \arccos r^T R^T v \), between the body mounted sensor and the inertially fixed sensor.
The experiment demonstrates that the minimum angular separation is \SI{14}{\degree} which satisfies the constraint of \( \theta = \SI{12}{\degree} \).
This further validates that the proposed control system exhibits the desired performance in the experimental setting as well. 
A video clip is available at \url{https://youtu.be/dsmAbwQram4}.

\section{Conclusions}\label{sec:conclusions}
We have developed a geometric adaptive control system which incorporates state inequality constraints on \(\SO\).
The presented control system is developed directly on \(\SO\) and it avoids singularities and ambiguities that are inherent to attitude parameterizations.
The attitude configuration error is augmented with a barrier function to avoid the constrained region, and an adaptive control law is proposed to cancel the effects of uncertainties. 
We show the stability of the  proposed control system through a rigorous mathematical analysis.
In addition, we have demonstrated the control system via numerical simulation and hardware experiments on a hexrotor UAV.
A novel feature of this control is that it is computed autonomously on-board the UAV.
This is in contrast to many state constrained attitude control systems which require an a priori attitude trajectory to be calculated. 
The presented method is simple, efficient and ideal for hardware implementation on embedded systems.

\appendix
\subsection{Proof of~\Cref{prop:config_error}}\label{proof:config_error}
To prove~\cref{item:prop_psi_psd}, we note that~\cref{eqn:A} is a positive definite function about \( R = R_d \)~\cite{bullo2004}.
The constraint angle is assumed \( \ang{0} \leq \theta \leq \ang{90} \) such that \( 0 \leq \cos \theta \).
The term \( r^T R^T v \) represents the cosine of the angle between the body fixed vector \( r \) and the inertial vector \( v \). 
It follows that
\begin{align*}
	0 \leq  \frac{\cos \theta -  r^T R^T v}{1 + \cos \theta} \leq 1 ,
\end{align*}
for all \( R \in \SO \). 
As a result, its negative logarithm is always positive and from~\cref{eqn:B}, \(1 < B\).
The error function \( \Psi = A B \) is composed of two positive terms and is therefore also positive definite.

Next, we consider~\cref{item:prop_era}.
The infinitesimal variation of a rotation matrix is defined as
\begin{align*}
    \delta R = \left. \diff{}{\epsilon} \right|_{\epsilon=0} R \exp{\epsilon \hat{\eta}} = R \hat{\eta} .
\end{align*}
Using this, the variation of~\cref{eqn:A} is taken with respect to \( R \) as
\begin{align*}
	\dirDiff{A}{R} &= \eta \cdot \frac{1}{2} \left( G R_d^T R - R^T R_d G\right)^\vee ,
\end{align*}
where we used~\cref{eqn:hat1} to achieve the simplified form.

A straightforward application of the chain and product rules of differentiation allows us to show~\cref{item:prop_erb} as
\begin{align*}
	\dirDiff{B}{R} &=  \eta \cdot \frac{ - \left( R^T v\right)^\wedge r}{\alpha \left(\cos \theta - r^T R^T v \right)} ,
\end{align*}
where the scalar triple product~\cref{eqn:STP} was used.

We show~\cref{item:prop_psi_quadratic} by computing the hessian of \( \Psi \), namely \( \Hess(\Psi) \),  using the chain rule as 
\begin{align*}
    \Hess &(\Psi)\cdot (\delta R,\delta R)= (\mathbf{D}_R\parenth{\mathbf{D}_R A \cdot \delta R}\cdot \delta R) B\\
    & + \parenth{\dirDiff{A}{R} }\parenth{\dirDiff{B}{R}} 
    + \parenth{\dirDiff{B}{R} }\parenth{\dirDiff{A}{R}}\\ &+(\mathbf{D}_R\parenth{\mathbf{D}_R B \cdot \delta R}\cdot \delta R) A.
\end{align*}
The first order derivative of \( A(R) \) and \( B(R) \) are given by~\cref{eqn:eRA,eqn:eRB}. 
The hessian of \( A(R) \) is computed as 
\begin{align*}
    \dirDiff{\parenth{\dirDiff{A}{R}}}{R} &= \eta \cdot \frac{1}{2}\parenth{G R_d^T \delta R - \delta R^T R_d G}^\vee  \\
    &= \eta \cdot \frac{1}{2} \parenth{G R_d^T R \hat \eta + \hat \eta R^T R_d G}^\vee\\
    &= \eta \cdot \frac{1}{2} \bracket{\parenth{\braces{\tr{R^T R_d G} I - R^T R_d G} \eta}^\wedge}^\vee  \\
    &= \eta \cdot \frac{1}{2} \parenth{\tr{R^T R_d G} I - R^T R_d G}\eta ,
\end{align*}
where we used the scalar triple product rule, from~\cref{eqn:STP}, to arrive at the final form.
The hessian of \( B(R) \) is computed as
\begin{align*}
    \dirDiff{\parenth{\dirDiff{B}{R}}}{R} =& \eta \cdot \left[\frac{\parenth{\delta R^T v}^\wedge r}{\alpha\parenth{r^T R^T v - \cos \theta}} \right. \\
    & \left. - \frac{\parenth{R^T v}^\wedge r \parenth{r^T \delta R^T v}}{\alpha \parenth{r^TR^Tv - \cos \theta}^2}   \right] .
\end{align*}
The term \( \parenth{\delta R^T v}^\wedge r\) is simplified as
\begin{align*}
    \parenth{\delta R^T v}^\wedge r &= - \hat{r} \delta R ^T v  \\
    &= - \hat{r} \parenth{R \hat \eta}^T v  \\
    &= \hat{r} \parenth{\hat{\eta} R^T } v  \\
    &= - \hat{r} \parenth{R^T v}^\wedge \eta ,
\end{align*}
where we utilized the hat map property from~\cref{eqn:cross_product}.
Similarly, the term \( r^T \delta R^T v \) is simplified as
\begin{align*}
    r^T \delta R^T v &= r^T \parenth{-\hat \eta R^T}v  \\
    &= r^T \parenth{R^T v}^\wedge \eta .
\end{align*}
The hessian of \( B(R) \) then becomes
\begin{align*}
    \dirDiff{\parenth{\dirDiff{B}{R}}}{R} =& \eta \cdot \left[ - \frac{ \hat r \parenth{R^T v}^\wedge}{\alpha \parenth{r^T R^T v - \cos \theta}} \right. \\ 
    & \left. - \frac{\parenth{R^T v}^\wedge r r^T \parenth{R^T v}^\wedge }{ \alpha\parenth{r^T R^T v - \cos \theta}^2}  \right] \eta .
\end{align*}
Using these terms, we evaluate \( \Hess{\Psi} \) at the desired attitude \( R = R_d \) as follows. Since $A=0$ and $\mathbf{D}_R A=0$ at $R=R_d$, 
\begin{align*}
    \Hess(\Psi)\cdot (\delta R,\delta R)\big|_{R=R_d} = \eta \cdot \frac{1}{2} B \parenth{\tr{G}I -  G} \eta , 
\end{align*}
which is positive definite since \( B > 1\) and \( \sum g_i > g_i\). 
The domain \( D \) is an open neighborhood of the desired attitude \( R_d \), and it excludes the undesired equilibrium points of \( A(R) \) and the infeasible regions defined by the constraints \( r^T R^T v_i \). 
Therefore, the only critical point of the error function $\Psi$ in the domain $D$ corresponds to the desired attitude $R=R_d$ with $e_R=0$ and $\Psi=0$. 
Therefore, in \( D \) the configuration error function is quadratic and the bounds in~\cref{item:prop_psi_quadratic} are valid according to~\cite[Proposition 6.30]{bullo2004}.

The proof of \cref{item:prop_era_upbound} is available in~\cite{LeeITCST13}.

\subsection{Proof of~\Cref{prop:error_dyn}}\label{proof:error_dyn}
From the kinematics~\cref{eqn:Rdot}, and noting that \( \dot{R}_d = 0 \) the time derivative of \( R_d^T R \) is given as
\begin{gather*}
	\diff{}{t} \parenth{R_d^T R} = R_d^T R \hat{e}_\Omega .
\end{gather*}
Applying this to the time derivative of~\cref{eqn:A} gives
\begin{gather*}
	\diff{}{t} (A) = -\frac{1}{2} \tr{G R_d^T R \hat{e}_\Omega} .
\end{gather*}
Applying~\cref{eqn:hat1} into this shows~\cref{eqn:A_dot}.
Next, the time derivative of the repulsive error function is given by
\begin{gather*}
	\diff{}{t} (B) = \frac{r^T \parenth{\hat{\Omega} R^T} v}{\alpha \parenth{r^T R^T v - \cos \theta}} .
\end{gather*}
Using the scalar triple product, given by~\cref{eqn:STP}, one can reduce this to~\cref{eqn:B_dot}.
The time derivative of the attractive attitude error vector, \( e_{R_A} \), is given by
\begin{gather*}
	\diff{}{t} ( e_{R_A}) = \frac{1}{2} \parenth{\hat{e}_\Omega R^T R_d G + (R^T R_d G)^T \hat{e}_\Omega}^\vee .
\end{gather*}
Using the hat map property given in~\cref{eqn:xAAx} this is further reduced to~\cref{eqn:eRA_dot,eqn:E}.

We take the time derivative of the repulsive attitude error vector, \( e_{R_B} \), as
\begin{gather*}
	\diff{}{t}( e_{R_B} )= a \Omega v^T R r - a R^T v \Omega^T r + b R^T \hat{v} R r ,
\end{gather*}
with \( a \in \R \) and \( b \in \R\) given by 
\begin{gather*}
	a = \bracket{\alpha \parenth{r^T R^T v - \cos \theta}}^{-1} , \;
	b = \frac{r^T \hat{\Omega} R^T v}{\alpha \parenth{r^T R^T v - \cos \theta}^2} .
\end{gather*}
Using the scalar triple product from~\cref{eqn:STP} as \( r \cdot \Omega \times \parenth{R^T v} = \parenth{R^T v} \cdot r \times \Omega \) gives~\cref{eqn:eRB_dot,eqn:F}.

We show the time derivative of the configuration error function as
\begin{gather*}
	\diff{}{t} (\Psi) = \dot{A} B + A \dot{B} .
\end{gather*}
A straightforward substitution of~\cref{eqn:A_dot,eqn:B_dot,eqn:A,eqn:B} into this and applying~\cref{eqn:eR} shows~\cref{eqn:psi_dot}.
We show~\cref{eqn:eW_dot} by rearranging~\cref{eqn:Wdot} as 
\begin{align*}
	\diff{}{t} e_\Omega = \dot{\Omega} = J^{-1} \parenth{u - \Omega \times J \Omega + W(R,\Omega) \Delta } .
\end{align*}

\subsection{Proof of~\Cref{prop:att_control}}\label{proof:att_control}
Consider the following Lyapunov function:
\begin{gather*}
	\mathcal{V} = \frac{1}{2} e_\Omega \cdot J e_\Omega + k_R \Psi(R,R_d) . 
\end{gather*}
From~\cref{item:prop_psi_psd} of~\Cref{prop:config_error}, \(\mathcal{V} \geq 0 \).
Using~\cref{eqn:eW_dot,eqn:psi_dot} with \( \Delta = 0 \), the time derivative of \( \mathcal{V} \) is given by
\begin{align*}
	\dot{\mathcal{V}} &= -k_\Omega \norm{e_\Omega}^2 . 
\end{align*}
Since \( \mathcal{V} \) is positive definite and \( \dot{\mathcal{V}} \) is negative semi-definite, the zero equilibrium point \( e_R, e_\Omega \) is stable in the sense of Lyapunov. 
This also implies \( \lim_{t\to\infty} \norm{e_\Omega} = 0 \) and \( \norm{e_R} \) is uniformly bounded, as the Lyapunov function is non-increasing. From \refeqn{eRA_dot} and \refeqn{eRB_dot}, $\lim_{t\to\infty} \dot e_R =0$. 
One can show that \( \norm{\ddot{e}_R} \) is bounded.
From Barbalat's Lemma, it follows \( \lim_{t\to\infty}\norm{\dot{e}_R} = 0 \)~\cite[Lemma 8.2]{khalil1996}. 
Therefore, the equilibrium is asymptotically stable. 
	
Furthermore, since \( \dot{\mathcal{V}} \leq 0 \) the Lyapunov function is uniformly bounded which implies 
\begin{align*}
	\Psi(R(t)) \leq \mathcal{V}(t) \leq \mathcal{V}(0) .
\end{align*}
In addition, the logarithmic term in~\cref{eqn:B} ensures \( \Psi(R) \to \infty \) as \( r^T R^T v \to \cos \theta \).
Therefore, the inequality constraint is always satisfied given that the desired equilibrium lies in the feasible set.
	
\subsection{Proof of~\Cref{prop:eR_dot_bound}}\label{proof:eR_dot_bound}

Consider the open neighborhood $D$ of $R=R_d$ defined in~\Cref{prop:config_error}.
The proof of the upper bound of \( A(R) \) is given in~\cite{LeeITCST13}.
The selected domain ensures that the configuration error function is bounded \( \Psi < \psi \).
This implies that that both \( A(R) \) and \( B(R) \) are bounded by constants \( c_A c_B < \psi < h_1\).
Furthermore, since \( \norm{B} > 1 \) this ensures that \( c_A, c_B < \psi\) and shows~\cref{eqn:AB_bound}.

Next, we show~\cref{eqn:E_bound,eqn:F_bound} using the Frobenius norm.
The Frobenius norm \( \norm{E}_F \) is given in~\cite{LeeITCST13} as
\begin{gather*}
	\norm{E}_F = \sqrt{\tr{E^T E}} = \frac{1}{2} \sqrt{\tr{G^2} + \tr{R^T R_d G}^2} .
\end{gather*}
Applying Rodrigues' formula and the Matlab symbolic toolbox, this is simplified to
\begin{gather*}
	\norm{E}^2_F \leq \frac{1}{4} \parenth{\tr{G^2} + \tr{G}^2} \leq \frac{1}{2} \tr{G}^2 ,
\end{gather*}
which shows~\cref{eqn:E_bound}, since \( \norm{E} \leq \norm{E}_F \).

To show~\cref{eqn:F_bound}, we apply the Frobenius norm \( \norm{F}_F \):
\begin{align*}
	\norm{F}_F =& \frac{1}{\alpha ^2 \parenth{r^T R^T v - \cos \theta}^2} \left[\tr{a^T a} - 2 \tr{a^T b} \right. \\
	&\left.+ 2 \tr{a^T c } + \tr{b^T b}  - 2 \tr{b^T c} + \tr{c^T c}\right] .
\end{align*}
where the terms \( a, b, \text{ and } c \) are given by
\begin{gather*}
	a = r^T R r I , \quad	b = R^T v r^T , \quad c = \frac{R^T \hat{v} R r v^T R \hat{r}}{r^T R^T v - \cos \theta}.
\end{gather*}
A straightforward computation of \( a^T a \) shows that
\begin{gather*}
	\tr{a^T a} = \parenth{v^T R r}^2 \tr{I} \leq 3 \beta^2 ,
\end{gather*}
where we used the fact that \( v^T R r = r^T R^T v < \beta \) from our given domain.
Similarly, one can show that \( \tr{a^T b} \) is equivalent to
\begin{gather*}
	\tr{a^T b} = v^T R r \tr{R^T v r^T} = \parenth{v^T R r}^2 \leq \beta^2 ,
\end{gather*} 
where we used the fact that \( \tr{x y^T} = x^T y \).
The product \( \tr{a^T c} \) is given by
\begin{gather*}
	\tr{a^T c} = \frac{v^T R r}{r^T R^T v - \cos \theta} \tr{\parenth{R^T v}^\vee \parenth{r v^T R} \hat{r} } ,
\end{gather*}
where we used the hat map property~\cref{eqn:RxR}.
One can show that \(\mathrm{tr}[a^T c] \leq 0 \) over the range \( -1 \leq v^T R r \leq \cos \theta \). 
Next, \( \tr{b^T b}\) is equivalent to
\begin{gather*}
	\tr{b^T b} = \tr{r v^T R R^T v r^T} = 1 ,
\end{gather*}
since \( r,v \in \Sph^2\).
Finally, \( \tr{c^T c} \) is reduced to
\begin{gather*}
	\tr{c^T c} = \tr{\hat{r} R^T v r^T \bracket{-I + R^T v v^T R} r v^T R \hat{r}} ,
\end{gather*}
where we used the fact that \( \hat{x}^2 = - \norm{x}^2 I + x x^T\).
Expanding and collecting like terms gives
\begin{gather*}
	\tr{c^T c } = \frac{1 - 2\parenth{v^T R r}^2 + \parenth{v^T R r}^4}{\parenth{r^T R^T v - \cos \theta}^2} . 
\end{gather*}
Using the given domain \( r^T R^T v \leq \beta \) gives the upper bound~\cref{eqn:F_bound}.
The bound on \( e_{R_A} \) is given in~\cref{eqn:psi_lower_bound} while \( e_{R_B} \) arises from the definition of the cross product \( \norm{a \times b} = \norm{a} \norm{b} \sin \theta \).
Finally, we can find the upper bound~\cref{eqn:eR_dot} as
\begin{gather*}
	\norm{\dot{e}_R} \leq \parenth{\norm{B} \norm{E} + 2 \norm{e_{R_A}} \norm{e_{R_B}} + \norm{A}\norm{F}} \norm{e_\Omega} \, .
\end{gather*}
Using~\cref{eqn:AB_bound} to \cref{eqn:eRB_bound} one can define \( H \) in terms of known values.

\subsection{Proof of~\Cref{prop:adaptive_control}}\label{proof:adaptive_control}
Consider the Lyapunov function \( \mathcal{V} \) given by
\begin{align*}
	\mathcal{V} = \frac{1}{2} e_\Omega \cdot J e_\Omega + k_R \Psi + c J e_\Omega \cdot e_R + \frac{1}{2 k_\Delta} e_\Delta \cdot e_\Delta , 
\end{align*}
over the domain \( D \), defined in~\cref{prop:config_error}. 
In this set, the properties of~\Cref{prop:config_error,prop:error_dyn} are satisfied.

From~\Cref{prop:config_error}, the configuration error function is locally quadratic and it is bounded in \( D \) by~\cref{eq:psi_bound}.
Using this, the Lyapunov function \( \mathcal{V} \) is bounded by
\begin{align*} 
	z^T W_1 z \leq \mathcal{V} \leq z^T W_2 z ,
\end{align*}
where \( e_\Delta = \Delta - \bar{\Delta} \), \( z = [\|e_R\|,\|e_\Omega\|,\|e_\Delta\|]^T\in\R^3 \) and the matrices \(W_1,W_2 \in \R^{3 \times 3}\) are given by
\begin{align*}
	W_1 & = \begin{bmatrix}
		k_R n_1 & -\frac{1}{2} c \lambda_M & 0 \\
		-\frac{1}{2} c \lambda_M & \frac{1}{2} \lambda_m & 0 \\
		0 & 0 & \frac{1}{2 k_\Delta}
	\end{bmatrix},\\
	W_2 & = \begin{bmatrix}
		k_R n_2 & \frac{1}{2} c \lambda_M & 0 \\
		\frac{1}{2} c \lambda_M & \frac{1}{2} \lambda_M & 0 \\
		0 & 0 & \frac{1}{2 k_\Delta}
	\end{bmatrix} .
\end{align*}
The time derivative of \( \mathcal{V}\) with the control input defined by~\cref{eqn:adaptive_control} is given as
\begin{align*}
	\dot{\mathcal{V}} =& - k_\Omega e_\Omega^T e_\Omega + \parenth{e_\Omega + c e_R}^T W e_\Delta - k_R c e_R^T e_R \nonumber\\
	&- k_\Omega c e_R^T e_\Omega + c J e_\Omega^T \dot{e}_R - \frac{1}{k_\Delta} e_\Delta^T \dot{\bar{\Delta}} , \label{eqn:vdot}
\end{align*}
where we used \( \dot{e}_\Delta = - \dot{\bar{\Delta}} \).
The terms linearly dependent on \( e_\Delta\) are combined with~\cref{eqn:delta_dot} to yield
\begin{align*}
	 e_\Delta^T \parenth{W^T \parenth{e_\Omega + c e_R} - \frac{1}{k_\Delta} \dot{\bar{\Delta}}} = 0 . 
\end{align*}
Using~\Cref{prop:eR_dot_bound}, an upper bound on \( \dot{\mathcal{V}} \) is written as
\begin{gather*}
	\dot{\mathcal{V}} \leq -\zeta^T M \zeta ,
\end{gather*}
where $\zeta=[\|e_R\|,\|e_\Omega\|]\in\R^2$, and the matrix \( M \in \R^{2 \times 2} \) is 
\begin{gather*}
	M = \begin{bmatrix}
		k_R c & \frac{k_\Omega c}{2} \\
		\frac{k_\Omega c}{2} & k_\Omega - c \lambda_M H
	\end{bmatrix} .
\end{gather*}

If \( c \) is chosen such that~\cref{eqn:c_bound} is satisfied then the matrices \( W_1, W_2 \) and \( M \) are positive definite.
This implies that $\mathcal{V}$ is positive definite and decrescent, and $\dot{\mathcal{V}}$ is negative semidefinite in the domain $D$. As such, the zero equilibrium is stable in the sense of Lyapunov, and all of the error variables are bounded. Furthermore, $\lim_{t\to\infty} \zeta=0$ according to the LaSalle-Yoshizawa theorem~\cite{khalil1996}. 

%
%

\bibliography{BibMaster,library}
\bibliographystyle{IEEEtran}

\end{document}